\documentclass[10pt,journal]{IEEEtran} \onecolumn 

\usepackage{amsmath,amsfonts}
\usepackage{algorithmic}
\usepackage{algorithm}
\usepackage{array}
\usepackage{textcomp}
\usepackage{stfloats}
\usepackage{verbatim}
\usepackage{graphicx}
\usepackage{cite}
\usepackage{hyperref}
\hypersetup{colorlinks=true}
\usepackage{booktabs,color}
\newcommand{\bfpara}[1]{\noindent{\bf #1. }}
\newcommand{\empara}[1]{\noindent{\em #1. }}
\begin{document}
%
\title{Democratizing Signal Processing and Machine Learning: Math Learning Equity for Elementary and Middle School Students}
%
\author{Namrata Vaswani \IEEEmembership{Fellow,IEEE}, Mohamed Y. Selim \IEEEmembership{Senior Member,IEEE}, Renee Serrell Gibert
\thanks{N. Vaswani is with the Electrical and Computer Engineering department at Iowa State University.  Email: namrata@iastate.edu. M. Selim is with the Electrical and Computer Engineering department at Iowa State University. R. Gibert is with the Minority Engineering Program at Purdue University.}%
}


\maketitle

{\em Submitted to IEEE Signal Processing Magazine.  
}

\begin{abstract}
Signal Processing (SP) and Machine Learning (ML) rely on good math and coding knowledge, in particular, linear algebra, probability, trigonometry, and complex numbers. A good grasp of these relies on scalar algebra learned in middle school. The ability to understand and use scalar algebra well, in turn, relies on a good foundation in basic arithmetic. Because of various systemic barriers, many students are not able to build a strong foundation in arithmetic in elementary school. This leads them to struggle with algebra and everything after that. Since math learning is cumulative, the gap between those without a strong early foundation and everyone else keeps increasing over the school years and becomes difficult to fill in college. In this article we discuss how SP faculty, students, and professionals can play an important role in starting, and participating in, university-run, or other, out-of-school math support programs to supplement students' learning. Two example programs run by the authors, CyMath at Iowa State and Algebra by 7th Grade (Ab7G) at Purdue, and one run by the Actuarial Foundation, are described.
We conclude with providing some simple zero-cost suggestions for public schools that, if adopted, could benefit a much larger number of students than what out-of-school programs can reach. 
\end{abstract}

\bfpara{About the authors} 
Vaswani is a Professor of Electrical and Computer Engineering (ECE), and Anderlik Professor of Engineering, at Iowa State University, where her teaching and research lie in statistical signal processing (SP) and machine learning (ML). She is a Fellow of IEEE (class of 2019) and of the AAAS (class of 2023). In the last 20 years, she has extensively published in various IEEE Transactions and has taught multiple Math for SP and ML courses. These include undergraduate Probability, an undergraduate course on ML from an SP perspective, graduate level Detection and Estimation Theory, and High-dimensional Probability for ML and SP. In 2020, she founded and directs and runs the CyMath (graduate student led) elementary and middle school mathematics tutoring and mentoring program. Selim is a networking and computer engineering educator and researcher and an Associate Teaching Professor at Iowa State since 2018. He has been a part of CyMath since 2024 and now leads its online tutoring component. Gibert is a Program Manager for the Minority Engineering Program at Purdue University and a former middle school mathematics teacher. She coordinates the Algebra by Seventh Grade (Ab7G) program at Purdue, which is a much larger scale math mentoring program for school students.%



\section*{Introduction}
This article explains why a good foundation in early math skills is critical for later success in signal processing (SP) and machine learning (ML). We discuss how early interventions can help reduce the learning gaps which keep growing over time; and hence promote equity in signal processing (SP) and machine learning (ML).
We first explain why early math is critical for SP and ML in Sec. \ref{whyearlymath}. This section can be appropriately modified, or skipped, and the rest of the article will still be relevant for most science, technology, engineering, and mathematics (STEM) fields. In Sec. \ref{barriers}, we describe (well-known and lesser known) systemic barriers that prevent some students from getting a good early math foundation. In Sec. \ref{soln1}, we describe how SP faculty, graduate students or professionals can play an important role in starting, and participating in, out-of-school math support programs to supplement students' learning.  We conclude in Sec. \ref{soln2} by providing some zero-cost suggestions that, if adopted by public schools, can help reduce learning gaps for a much larger number of students than what out-of-school programs can reach. 
These are based on our perspective as college educators teaching math and coding based engineering classes and who have also been involved in directing or running out-of-school math support programs for young students.


{\em The first section below can be appropriately modified or ignored and the rest of the article will still be relevant for most science, technology, engineering, and mathematics (STEM) fields.
}

\section{Why Early Math Education matters for Signal Processing (SP) and Machine Learning (ML)} \label{whyearlymath} 
%
Understanding modern signal processing (SP) and machine learning (ML) -- basic SP, statistical SP, digital SP, estimation and detection theory, or ML algorithms -- requires a good grasp of linear algebra (vector and matrix algebra),  probability and statistics, and complex numbers, and an ability to code algorithms based on these mathematics (math) concepts in MATLAB or Python. The latter is essential to understand, develop, and evaluate novel applications of SP and ML algorithms. However, one cannot learn linear algebra, or even learn to code it in, without a strong knowledge of basic scalar algebra that is typically taught in middle school, and sometimes in high school. Algebra also forms the foundation of the ability to understand abstract probability concepts such as probability mass function or probability density function. Algebra, in turn, cannot easily make sense to a student without fluency in elementary school arithmetic (add, subtract, multiply divide; negative numbers; fractions, decimals, and four operations with these; exponents). As an example, the ability to solve for $x,y$ from two equations in $x,y$: $ 0.2(3x+5) - 9x + 22.1= 1/3$ and $x-y=0.29$ requires arithmetic fluency. Solving for two variables from two linear equations forms the basis of understanding linear algebra, e.g., how to solve for an unknown vector $\bf{x}$ from $\bf{A} \bf{x} = \bf{b}$. This, in turn, is the basis for understanding Least Squares estimation which is used for parameter learning for Linear Regression. This understanding is essential to correctly code in the least squares estimation algorithm, or to correctly use software tools (e.g., the backslash operator in MATLAB or operators from the NumPy package in Python) that provide such a solution.
In fact, in the ChatGPT era, where code snippets and code syntax can be automatically generated, math will become even more important because one needs the skills to verify whether the auto-generated code is indeed mathematically correct.

As an example, implementing least squares estimation to solve for $\bf{x}$ from $\bf{b}:= \bf{A} \bf{x} + \bf{w}$, where $\bf{w}$ is modeling error and $\bf{A}$ is a tall matrix, requires computing $\bf{A}^\dag \bf{b} := (\bf{A}^\top \bf{A})^{-1} \bf{A}^\top \bf{b} $. This can be coded in various ways in Python, but some are better than others in certain settings. For instance, the pinv function from some of the Python packages does not always work as desired.  Moreover, the condition number of $\bf{A}$ governs solution accuracy, and when it is ill-conditioned, none of the approaches give an accurate answer.   ``Regularization'' approaches exist to deal with ill-conditioning. Realizing  all this requires a good understanding of linear algebra.

\begin{figure*}
  \centering
  \includegraphics[height=3in]{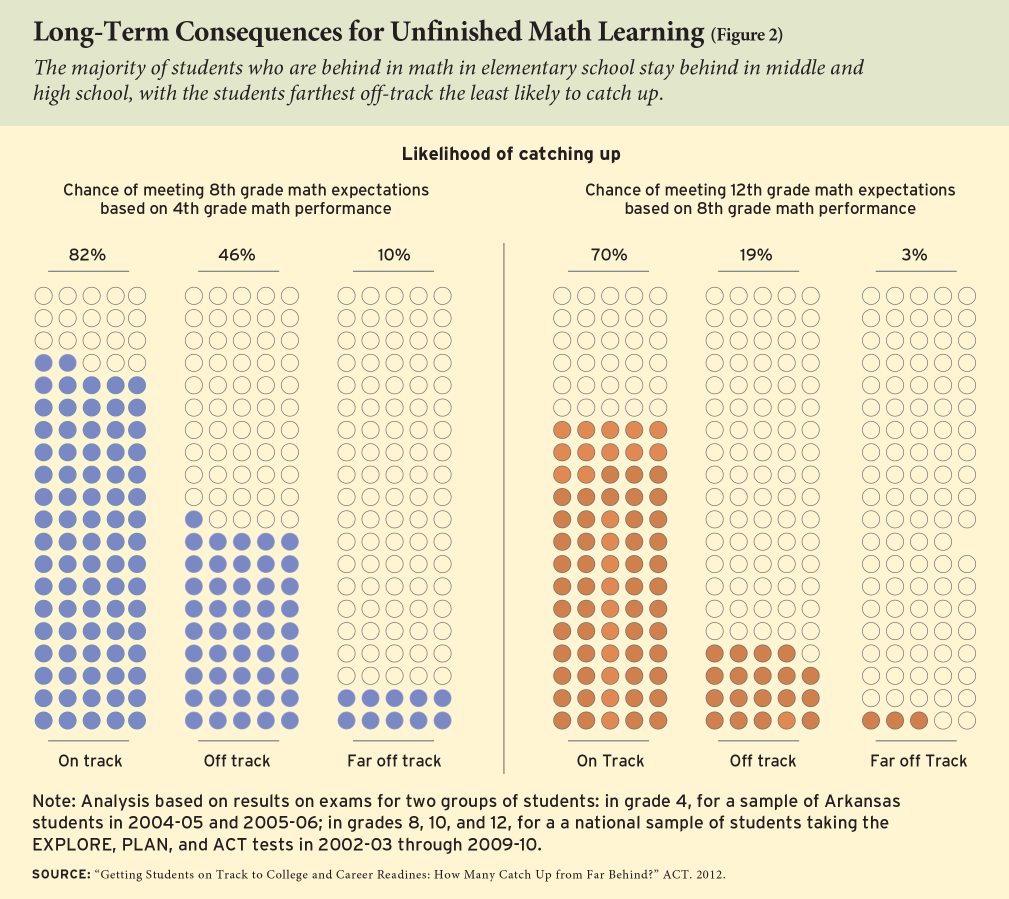}
    \vspace{-0.1in}
  \caption{\sl{Figure indicating that math learning is cumulative. This is taken from \url{https://www.educationnext.org/grade-level-expectations-trap-how-lockstep-math-lessons-leave-students-behind/}
which uses data from \cite{dougherty2012getting}. 
  }}
    \label{fig1}
  \vspace{-0.1in}
\end{figure*}

In summary, SP and ML rely on good college-level math knowledge. Since math learning is cumulative, the gap between those who did not get an opportunity to learn arithmetic, and hence scalar algebra, well and everyone else keeps increasing over the school years \cite{dougherty2012getting} and becomes difficult to fill by the time they reach college. We provide a visual in Fig. \ref{fig1} that shows this. 
One way to make SP learning more accessible is to fill the math learning gap while it is still small by providing extra support early. 



%

\section{Systemic Barriers} \label{barriers} 
Due to various systemic barriers, many students do not get a good early foundation in mathematics and consequently, opportunities for economic mobility via higher education are blocked to them \cite{moses2002radical}.
Well-known barriers include inadequate resources for public education, family income and education levels, non-uniform teacher quality, COVID-19 learning losses, unconscious bias,  language barriers and a culture of low expectations \cite{flores2007examining}\footnote{``Specifically, data show that African American, Latino, and low-income students are less likely to have access to experienced and qualified teachers, more likely to face low expectations, and less likely to receive equitable per student funding''.},\cite{goss2016widening}\footnote{``Learning gaps between Australian students of different backgrounds are alarmingly wide and grow wider as students move through school. Importantly, the learning gaps grow much larger after Year 3. 
Bright kids in disadvantaged schools show the biggest losses, making two-and-a-half years less progress than students with similar capabilities in more advantaged schools.''
}.
In addition to the above, the authors of this article have also noticed the following lesser-known barriers, which should be easier to remove since they are due to design changes that were introduced to improve the elementary school experience. These include:
(i) insufficient focus on math practice in elementary school (both during school hours or as homework);
(ii) little information provided to families that a good grasp of elementary school math is essential for later success, and how to help students get this grasp.
(iii) Most parents are not made aware that results of tests taken by students in fourth or fifth grade are used to decide which student gets to join advanced math tracks in fifth or sixth grade and on wards. Moreover in some school systems this track has been removed and this can be even worse for the high-performing students in these schools. 
%
%
As noted earlier, good early math skills are needed for learning algebra well.  Access to algebra has been referred to as the greatest civil rights issue of our time as it allows students access to economic ladders of opportunities  \cite{moses2002radical}. 

\section{A Partial Solution: Out-of-school Math Tutoring or Mentoring Programs} \label{soln1}
We describe three out-of-school university run math programs here -- \href{https://cymath.iastate.edu}{CyMath}, \href{https://www.purdue.edu/mep/Pre-College-Programs/ab7g\%20index.html}{Algebra by 7th Grade (Ab7G)}, and  \href{https://actuarialfoundation.org/programs-resources/math-motivators/}{Math Motivators}.
\href{https://cymath.iastate.edu}{CyMath} is Iowa State University's math tutoring program that was started by Vaswani in Fall 2020 and that primarily relies on volunteer graduate students or faculty as the tutors. Selim  leads CyMath's online tutoring component and has been with CyMath since early 2024. CyMath has been expanded during 2023-24 with help from Ab7G and Gibert. It is a small program that provides intensive (twice-a-week in-person or hybrid) math support along with at-home math resources. \href{https://www.purdue.edu/mep/Pre-College-Programs/ab7g\%20index.html}{Ab7G} is an elementary and middle school math mentoring program that was started at Purdue in 2017 by their Minority Engineering Program (MEP) and is coordinated by the third author, Renee Gibert. Ab7G provides monthly math mentoring followed by STEM exposure labs.  Most of its mentors are paid undergraduate students.
Both provide free regular math support, along with at-home learning resources, to students in grades 3-8 while prioritizing those from backgrounds  under-represented in Engineering and Computing. Math Motivators, run by the Actuarial Foundation, has the same goals but a slightly different approach and often runs in-school tutoring sessions.

We should emphasise that ``Cy'' (short for Cyclones) is ISU's mascot and hence the name CyMath was adopted. The CyMath program at ISU with website \url{https://cymath.iastate.edu} should not be confused with a completely unrelated high school math help company and website called cymath.com.

\subsection{CyMath: run by Electrical and Computer Engineering faculty at Iowa State University} \href{https://cymath.iastate.edu}{CyMath} 
was founded in Fall 2020 by Vaswani, with help from an Education colleague, Jackson, to help support math learning for students from disadvantaged backgrounds during COVID-19.
It started in online mode for one Des Moines, Iowa school in 2020 and ran that way until 2022, when interest in online-only tutoring fell. 
In 2023 Fall, with help and advice from Ab7G, Vaswani re-started CyMath as an in-person program for the local Ames Iowa school district.  
 CyMath begins tutoring students in third or fourth grade and follows them through the school years, allowing tutors to fill in the learning gaps while they are still small. 
Its long-term goal is to increase the number of youth prepared to pursue, and thrive in, Engineering, Computing or other STEM majors in college.

Most CyMath tutors are volunteer graduate students or faculty from Electrical and Computer Engineering, Mathematics, and Statistics, with some also from other Engineering departments. These are supported by a few paid undergraduate student tutors and assistants, most of whom are  Education majors, while some are STEM undergraduate students.  These tutors help run the program, help students (children) stay focused, and serve as tutors for certain students. Regular text messages and phone calls, along with email, are used to motivate and encourage parents to encourage their students to also do at-home math practice.
Borrowing from Ab7G, CyMath buys accounts for McGrawHill's adaptive math learning software ALEKS (Assessment and Learning in Knowledge Spaces) \url{www.aleks.com}. Tutors use this as a base for tutoring. Students are encouraged to regularly work on math on ALEKS at home. Roughly 20-30 minutes a day for about five days a week is recommended. For students who do not have reliable internet access, tutors try to send workbooks and worksheets (printed from websites such as \url{k5learning.com}) home. CyMath does not insist on a particular tutoring curriculum, or use of ALEKS, and so far has provided limited training to tutors. It relies on the assumption that tutors know their math (have graduate level training), have a passion for it, and are responsible adults (tutors are carefully chosen and observed).

A large fraction of CyMath tutors are international graduate students or immigrant faculty, including the program leads Vaswani and Selim. Since the language of math is universal, these tutors have been as effective as domestic graduate students. In fact, they also bring an international perspective to math learning and expose students to world cultures and math teaching approaches, that is often completely missing in elementary schools.
At the same time, the tutoring helps the tutors become better college educators (teaching assistants or instructors), and provides them an understanding of the modern US math education system. This enables them to pace their own college classes better.
%
Finally, when possible, CyMath attempts to inform parents and the students about the dates for state-wide testing and also what adaptive testing means (need to answer the early easy problems well to get the harder problems which are also worth more points).  Lastly, similar to Ab7G, CyMath runs engineering fun sessions roughly twice a semester.

%

In the online mode (2020-2022), tutoring was provided online to students in groups of two or three and once a week. CyMath supported 14 and 35 students in the two years. During 2023-24, tutoring was fully in-person and one-on-one. Each student got one session a week in Fall 2023, and two a week in Spring 2024. Tutoring was followed by a 20-minute play time with LEGOs or on the school playground. Tutoring ran after-school at the school itself and busing was provided by the district to take students back home, thus reducing demands on parents' time. CyMath benefited significantly from the school principal's support. During Summer 2024, Zoom sessions, and some in-person ones at Iowa State, were run. We taught a total of about 25 students for some part of the school year. New students were added as soon as more tutors became available or when some students left (September, January, March, May). For the 2024-25 school year, CyMath has transitioned to hybrid mode tutoring. Each week, one session is in-person and the second one is online for students with internet access at home. For those without, the second session is also in-person. 
Hybrid tutoring helps with reaching parents more easily, and maintaining tutor continuity after the tutor graduates. Finally, it is helping CyMath expand more easily to schools that are farther away. 

%

For the earliest group of six students who started in Fall 2023 when the local school district CyMath started, we observed the following in test score data on standardized tests conducted by the school -- in 2023-2024, the school used Measures of Academic Progress (MAP), while in 2024-25, they are using the iReady tests. One student has transitioned from being at a 40-th percentile in September 2023 to being at the 85th percentile in May 2024 and at the 92-nd percentile in January 2025. A second student from this cohort has been showing more modest but steady
gains and has gone from 20-th to 38-th to now 60-th percentile. One newer student who started in Fall 2024
has also gone up from 34-th to 66-th percentile, while some of the other newer students are going up by 5-10
percentile points. However, since percentile measures are noisy (the score gaps are very small), only large
and consistent changes are meaningful.

Lessons learned from CyMath include: (1)  Tutor-student need to connect well for tutoring to work; graduate students will often tutor for 2-3 years and hence their help can be quite useful. Use of hybrid-mode tutoring helps with this continuity. (2) A combination of mentors who are  Elementary Education majors (future school teachers) and Engineering/Math majors is an excellent way to run the sessions. Education majors are excellent at helping students who are really struggling with math or emotionally or both. STEM tutors bring a math learning approach that they see as essential for also succeeding in later math courses; this emphasises more math practice, working through worksheets, and not always converting each math exercise into a game. In fact, both types of tutors can influence each other in positive ways. 
(3) Encouraging students to complete their math work first, followed by the promise of game time or playground time or both, works better than vice-versa.

\subsection{Algebra by 7th Grade (Ab7G): run by Purdue University School of Engineering}
 \href{https://www.purdue.edu/mep/Pre-College-Programs/ab7g\%20index.html}{Ab7G} was established at Purdue in 2017. Chevron and Duke Energy seeded the funding for Ab7G with a generous grant. The mission of Ab7G is to increase the number of under-represented 7th grade students in the United States that are academically prepared to take algebra. The program is provided at no cost to participants.
Under the leadership of Director Virginia Booth Womack of the Minority Engineering Program (MEP) in 2017, forty-four students in grades three through seven from local school districts near Purdue's West Lafayette campus participated as the inaugural cohort. A partnership was created with Lafayette School Corporation (LSC) to help recruit students for the program. LSC principals and math coordinators were very helpful in distributing flyers and allowing MEP to attend math and science nights to recruit students. At the conclusion of each academic year, MEP provides the district office with comprehensive reports on students' final performance and attendance for the sessions.
Forty-four students attended in-person sessions twice per month for a total of sixteen sessions during the 2017-2018 academic calendar. By 2024, Ab7G had 208 participants, including some virtual-only participants from outside Indiana and outside the US also.

The exigencies of the COVID-19 pandemic forced the program to move to a virtual platform. Currently, the program operates in a hybrid format, offering two in-person sessions each semester, accompanied by an asynchronous virtual option during the in-person sessions. Additionally, four virtual sessions are provided each semester. The Saturday sessions are structured for students to work on McGraw Hill's ALEKS  (Assessment and Learning in Knowledge Spaces) adaptive mathematics learning  software, mentorship with an undergraduate student, in-person sessions lunch, and a fun-lab (STEM activity).
Ab7G participants are encouraged to work at least fifteen minutes per day most days a week on ALEKS outside of the Saturday program and attempt to master three topics per week. At the beginning of each session, students are recognized for reaching their time and topic goals and honorable mentions are given to students that come close.

The program is built upon three foundational pillars: student self-efficacy, mentorship, and parental engagement. Students enhance their mathematical proficiency through both online and in-person instruction, collaborative activities, practical projects, and opportunities to interact with qualified STEM experts. 
Parent workshops involve exploring mathematical principles, navigating the online math tool, completing practical projects specifically created to enhance their child's learning, developing awareness for adolescent mental health, and creating growth mindsets.
Also, critical in the growth of the students in the fun lab curriculum.  Ab7G's fun lab curriculum covered topics on sustainable energy and building Lego windmills, planning a garden with measurement, Scratch Cat coding, Microbit activities, On Shape, 3-D printing, and financial literacy.
Gibert et al. \cite{gibert2024} discovered that utilizing a financial literacy curriculum that infused students' cultural experiences in a virtual platform changed their perceptions of money by understanding economic equity.



Lessons learned over the seven years of the Algebra by 7th Grade (Ab7G) program are as follows: (1) establishing strong partnerships with schools is crucial for the recruitment and retention of participants; (2) ensuring access to technology and high-speed internet is fundamental to the program's effectiveness and to the academic success of the students; (3) providing translations for documents is useful; and (4) designing the program in a way that encourages the retention of mentors each year is essential for continuity. (4) Womack et al \cite{booth2022program} reported significant growth in mathematics proficiency among students who completed a full grade level in ALEKS each year, n=20. This growth demonstrated a direct correlation with achieving proficiency or higher on Indiana's Learning Evaluation and Assessment Readiness Network (ILEARN) test which is administered to Indiana students in grades 3-8 at the end of each school year.%



\subsection{Math Motivators: run by the Actuarial Foundation}
\href{https://actuarialfoundation.org/programs-resources/math-motivators/}{Math Motivators} is a math tutoring program run by the Actuarial Foundation and provides in and out of school math tutoring support to students from disadvantaged backgrounds in various schools around the US.  It is supported by the Society of Actuaries (a world-wide organization of actuaries), the American Academy of Actuaries, and various corporate sponsors.
Approximately 75\% of its participating students ares eligible for free or reduced school lunch, and 70\% identifying as Black, African American, Latino or of Spanish origin.
Math Motivators offers free in-person math tutoring to under-served students in grades 3-12 who need and want tutoring but otherwise cannot afford it. Two students of similar ability are paired with one tutor, who works with the students to help them become proficient in math. Some of their programs offer tutoring during the school day while others offer tutoring after-school. MM currently runs in eight cities in the USA.

It would be interesting to study the impact if IEEE, or a society within IEEE, partners with Math Motivators, or starts its own Math and Coding Support program. IEEE likely has a much larger membership base than that of actuarial societies.





\section{Lessons learned and Conclusions}\label{soln2}
One goal of writing this article is to invite SP readers from around the world to start, participate in, and write about math support programs for their local schools. 
Our hope is also that this article encourages universities/colleges and professional organizations to support their faculty  or their membership interested in starting such programs. Of course these could also be started by community or religious organizations that reach out to colleges for tutors.
%
Out-of-school programs cannot reach all the students that public schools reach. Even if there were infinite resources available to teach and transport all children, many will not want to, or will not be able to, join one due to other commitments. 
A second goal of this article is to share  the lessons we have learned as college math and engineering educators who have also been running the school math support programs. 
While school math curriculums differ significantly across the world, and even within the same country, the math background needed to succeed in college signal processing (SP) courses is roughly the same everywhere; see  Sec. \ref{whyearlymath}. Furthermore the skills needed to succeed as an SP professional are also similar worldwide.

%
%


\begin{itemize}
\item  \empara{Math Homework and Sufficient In-School Math Practice}
Individual math practice is essential for math fluency. Thus, sending reasonable amounts of math practice work home as homework can have significant benefits, especially when teachers have the time to check the students' work carefully and discuss it with them so that student learning can be corrected early and easily.  College courses, and often high school courses, do have significant homework and this would also teach students an excellent skill needed for college success. The literature on elementary school homework is mixed. Some studies, such as  \cite{homework_latin_america} (based on Latin American schools) or \cite{homework_china} (based on China's schools) argue that reasonable amounts of math homework in primary (elementary) school does, in fact, improve student achievement. The latter also shows that too much homework does not provide extra benefit. 
On the other hand, other  studies such as \cite{jerrim2020association} conclude that homework does not help in elementary school. However this study followed students for a short period of time (less than two years). This is the case for a lot of elementary education literature since it is not easy to do long-term studies in countries where school systems, and hence student data, are decentralized. Such studies miss the long-term (elementary school teaching policies influencing middle and high school and college outcomes) effects of their proposed changes. 

Besides homework, providing enough in-school individual math practice can also help with student learning in younger grades. This would be easier for teachers to incorporate if they were not encouraged convert every math lesson into a game. Planning for fun math games takes a lot of teacher time and uses up a lot of classroom time. It likely takes away the class time that could otherwise be spent by the students on math practice; and by the teachers on checking the students' work and correcting errors early.



\item \empara{Parent Awareness}
Many parents are not aware that a good early foundation in math is important, or how to help their child build such a foundation. 
This knowledge should not be the privilege of children or relatives of STEM professionals. It should be provided to everyone by the schools. 
Moreover, it will automatically get conveyed if math is emphasised, homework is sent, and if free workbooks, printed worksheets, and information about good online resources such as Khan academy or k5learning, are shared with families. Some information is summarized here \href{https://cymath.iastate.edu/math-for-all/}{Math-info}.



Moreover, in many elementary schools, there are few tests, and little feedback is provided to parents based on the test results. Parents or students are often not even informed about the test dates even though the same tests' results are used later to place students into different math tracks.
However, if a parent asks, they are provided all of this information. Consequently, we end up disadvantaging a much larger number of children who may have otherwise come better prepared for these tests. Good test taking skills are also quite important for college success.

\item \empara{Use of summer to make-up for lost learning}  Despite all efforts, some students may require additional support to stay on track.
The summer before middle school can be used efficiently to help students catch up on arithmetic and pre-algebra skills for example.
Summer is also the time many undergraduate and graduate students and faculty are more available and schools can tap into these resources by partnering with nearly colleges or universities.



\end{itemize}





\bibliographystyle{IEEEtran}
\bibliography{../CyMathProposals/citations,../../bib/tipnewpfmt_kfcsfullpap}

\end{document}